\documentclass[3pt]{article}
\renewcommand{\baselinestretch}{1}

\usepackage{graphicx}\usepackage{epstopdf}\usepackage{color}
\usepackage{amsmath}\usepackage{amsfonts}\usepackage{amssymb}
 \usepackage{amsthm}

\numberwithin{equation}{section}

\begin{document}

\title{An adaptive algorithm based on the shifted inverse iteration for the Steklov eigenvalue problem }

\author{ {Hai Bi, Hao Li, Yidu Yang} \\\\
{\small School of Mathematics and Computer Science, }\\{\small
Guizhou Normal University,  Guiyang,  $550001$,  China}\\{\small
 bihaimath@gznu.edu.cn,lihao365@126.com,  ydyang@gznu.edu.cn}}
\date{~}
\pagestyle{plain} \textwidth 145mm \textheight 215mm \topmargin 0pt
\maketitle

\indent{\bf\small Abstract~:~} {\small This paper proposes and
analyzes an a posteriori error estimator for the finite element
multi-scale discretization approximation of the Steklov eigenvalue
problem. Based on the a posteriori error estimates, an adaptive
algorithm of shifted inverse iteration type is designed. Finally,
numerical experiments comparing the performances of three kinds of
different adaptive algorithms are provided, which illustrate the
efficiency of the adaptive algorithm proposed here.
}\\
\indent{\bf\small Keywords~:~\scriptsize} {\small Steklov eigenvalue
problem, finite element, multi-scale discretization, adaptive
algorithm, a posteriori error estimate.} 



\section{Introduction}
\indent In recent years, numerical methods for Steklov eigenvalue
problems have attracted more and more scholars' attention (see,
e.g.,\cite{aa1,andreev,armentano2,armentano3,bi1,cao,gm,huang,liq2,li,lx,aa2,xie}).
It is well known that in the numerical approximation of partial
differential equations, the adaptive procedures based on a
posteriori error estimates, due to the less computational cost and
time, are the mainstream direction and have gained an enormous
importance. The aim of this paper is to propose and analyze an a
posteriori error estimator for the finite element multi-scale
discretization approximation of the Steklov eigenvalue problem,
based on which an adaptive algorithm is designed.

\indent As for eigenvalue problems, till now, there are basically
three ways to design adaptive algorithms as follows in which the a
posterior error estimators are more or less the same but the
equations solved in each iteration are different: I. Solve the
original eigenvalue problem at each iteration. The convergence and
optimality of this adaptive procedure has been studied in
\cite{dai2,dai1,gm}. II. Inverse iteration type.
\cite{carstensen1,dr,mm,rannacher,rs,pssg} have studied and obtained
the convergence of this method. III. Shifted inverse iteration type
(see \cite{han,li2,yang3}).
 This paper studies the third type of adaptive method for the Steklov eigenvalue problem, and the special features are:\\
(1)~As for the Steklov eigenvalue problem, so far, it has been discussed the first type of adaptive method of combining the a posteriori error estimates and adaptivity (see \cite{gm} or Algorithm 4.1 in this paper). As far as our information goes, there has not been any report on the other two kinds of adaptive methods.
This paper designs the third type of adaptive algorithm based on the a posterior error estimates (see Algorithm 4.3 in this paper). Here we propose an a posteriori error estimator of residual type and give not only the global upper bound but also the local lower bound of the error which is important
for the adaptive procedure. \\
(2)~ \cite{li2} established and analyzed the a posteriori error estimates of the multi-scale discretization scheme for second order self-adjoint elliptic eigenvalue problems with homogenous Dirichlet boundary value condition by means of the a posteriori error estimates of the associated boundary value problem and left the local lower bound of the error unproved, while this paper studies the a posteriori error estimates of  the multi-scale discretization scheme by using the a posteriori error estimates of finite element eigenfunctions directly, and obtains the local lower bound of the error. \\
(3)~Numerical experiments comparing the performances of three kinds of different adaptive methods mentioned above are provided. It can be seen from the numerical results that the adaptive algorithm of shifted inverse iteration type has advantages over the other two kinds. More precisely, comparing with the first type, to achieve the same accurate approximation, our method uses less computational time; and our adaptive algorithm can be used to seek efficiently approximations of any eigenpair of the Steklov eigenvalue problem, however, the algorithm of the second type (see Algorithm 4.2 in this paper) is only suitable to the smallest eigenvalue.

\indent The rest of this paper is organized as follows. In section 2,
some preliminaries needed in this paper are introduced. In section 3, a multi-scale  discretization
scheme is presented, and its a priori and a posteriori error estimates are given and analyzed, respectively.
In section 4, three kinds of adaptive algorithms for the Steklov eigenvalue problem are presented, and finally numerical experiments are provided which illustrate the advantages of our algorithm.

\setcounter{section}{1}\setcounter{equation}{0}

\section{Preliminaries}
\indent Let $H^{t}(\Omega)$ and $H^{t}(\partial \Omega)$ denote Sobolev spaces on $\Omega$ and $\partial \Omega$ with real order $t$, respectively. The norm in $H^{t}(\Omega)$ and $H^{t}(\partial \Omega)$ are denoted by $\|\cdot\|_{t}$ and $\|\cdot\|_{t,\partial\Omega}$, respectively. $H^{0}(\partial\Omega)=L_{2}(\partial\Omega)$.\\
In this paper, we will write $a\lesssim b$ to indicate that $a\leqslant C b$ with $C>0$ being a constant depending on the data of the problem but independent of meshes generated by the adaptive algorithm.\\
\indent We consider the following Steklov eigenvalue problem
\begin{eqnarray}\label{s2.1}
 -\Delta u+u=0~~in~ \Omega,~~~~
\frac{\partial u}{\partial n}=\lambda u~~ on~ \partial\Omega,
\end{eqnarray}
where $\Omega \subset \mathbb{R}^{2}$ is a polygonal domain with
$\theta$ being the largest inner angle of $\Omega$ and
$\frac{\partial u}{\partial n}$ is the outward normal derivative.\\
\indent The weak form of (\ref{s2.1}) is given by: find $\lambda\in
\mathbb{R}$, $u \in H^{1}(\Omega)$ with $\|u\|_{1}=1$, such that
\begin{eqnarray}\label{s2.2}
a(u,v)=\lambda b(u,v),~~~\forall v\in H^{1}(\Omega),
\end{eqnarray}
where
\begin{eqnarray*}
a(u,v)&=&\int\limits_{\Omega}\nabla u\cdot \nabla v+uv
dx,~b(u,v)=\int \limits_{\partial\Omega}uvds ,\\
\|u\|_{b}&=&b(u,u)^{\frac{1}{2}}=\|u\|_{0,\partial\Omega}.
\end{eqnarray*}
It is easy to know that $a(\cdot,\cdot)$ is a symmetric, continuous
and $H^{1}(\Omega)$-elliptic bilinear form
on $H^{1}(\Omega)\times H^{1}(\Omega)$. So, we use $a(\cdot,\cdot)$ and
$\|\cdot\|_{a}=\sqrt{a(\cdot,\cdot)}=\|\cdot\|_{1}$ as the inner
product and norm on $H^{1}(\Omega)$, respectively.\\
\indent Let $\{\pi_{h}\}$ be a family of regular triangulations of $\Omega$ with the mesh
diameter $h$, and $V_{h}\subset C(\overline{\Omega})$ be a space of
piecewise linear polynomials defined on $\pi_{h}$.\\
 \indent The conforming finite element approximation of
(\ref{s2.2}) is: find $\lambda_{h}\in \mathbb{R}$, $u_{h} \in V_{h}$ with
$\|u_{h}\|_{a}=1$, such that
\begin{eqnarray}\label{s2.3}
a(u_{h},v)=\lambda_{h} b(u_{h},v),~~~\forall v\in V_{h}.
\end{eqnarray}

\indent Consider the following source problem (\ref{s2.4}) associated with
(\ref{s2.2}) and the approximate source problem (\ref{s2.5}) associated with (\ref{s2.3}), respectively.\\
Find $w \in H^{1}(\Omega)$, such that
\begin{eqnarray}\label{s2.4}
a(w,v)=b(f,v),~~~\forall v \in H^{1}(\Omega).
\end{eqnarray}
Find $w_{h} \in V_{h}$, such that
\begin{eqnarray}\label{s2.5}
a(w_{h},v)=b(f,v),~~~\forall v \in V_{h}.
\end{eqnarray}

From \cite{dauge} we know that the following regularity estimates hold for (\ref{s2.4}).\\
\indent{\bf Lemma 2.1.} If $f\in L_{2}(\partial\Omega)$, then there
exists a unique solution $w\in H^{1+\frac{r}{2}}(\Omega)$ to
(\ref{s2.4}), and
\begin{eqnarray}\label{s2.6}
\|w\|_{1+\frac{r}{2}}\lesssim \|f\|_{0,\partial\Omega};
\end{eqnarray}
if $f\in H^{\frac{1}{2}}(\partial\Omega)$, then there exists a unique solution $w\in
H^{1+r}(\Omega)$ to (\ref{s2.4}), and
\begin{eqnarray}\label{s2.7}
\|w\|_{1+r}\lesssim \|f\|_{\frac{1}{2},\partial\Omega},
\end{eqnarray}
where $r=1$ if $\Omega$ is convex, and $r<\frac{\pi}{\theta}$ which
can be arbitrarily close to $\frac{\pi}{\theta}$ when $\Omega$ is
concave.\\

\indent Then, thanks to Lemma 2.1, from (\ref{s2.4}) we can define
the operator $A: L_{2}(\partial\Omega)\rightarrow H^{1}(\Omega)$ by
\begin{eqnarray}\label{s2.9}
a(Af,v)=b(f,v),~~~\forall v \in H^{1}(\Omega).
\end{eqnarray}
\noindent Similarly, from (\ref{s2.5}) we define the operator
$A_{h}: L_{2}(\partial\Omega)\rightarrow V_{h}$ by
\begin{eqnarray}\label{s2.10}
a(A_{h}f,v)=b(f,v),~~~\forall v \in V_{h}.
\end{eqnarray}

It is obvious that $A: H^{1}(\Omega)\rightarrow H^{1}(\Omega)$ is a
self-adjoint operator. In fact, for any $u,v\in H^{1}(\Omega)$,
$$a(Au,v)=b(u,v)=b(v,u)=a(Av,u)=a(u,Av).$$
Analogously, $A_{h}$ is
also a self-adjoint operator. Observe that $Af$ and $A_{h}f$ are the
exact solution and the finite element solution of (\ref{s2.4}),
respectively, and $$a(Af-A_{h}f,v)=0, \forall v\in V_{h}\subset
H^{1}(\Omega).$$
Define the Ritz-Galerkin projection operator $P_{h}:
H^{1}(\Omega)\rightarrow V_{h}$ by
\begin{eqnarray*}
a(u-P_{h}u,v)=0,~~~\forall u \in H^{1}(\Omega),\forall v \in V_{h}.
\end{eqnarray*}
Thus, for any $f \in H^{1}(\Omega)$,
\begin{eqnarray*}
a(A_{h}f-P_{h}(Af),v)=a(A_{h}f-Af+Af-P_{h}(Af),v)=0,~~~ \forall v\in V_{h},
\end{eqnarray*}
therefore, $A_{h}f=P_{h}Af, \forall f \in H^{1}(\Omega)$, then $A_{h}=P_{h}A$.\\
From Lemma 2.1 and the interpolation error estimate, we deduce
\begin{eqnarray}\label{s2.132}
\|A_{h}-A\|_{a}&=&\sup\limits_{g\in
H^{1}(\Omega)}\frac{\|(A_{h}-A)g\|_{a}}{\|g\|_{a}}\nonumber\\
&=& \sup\limits_{g\in
H^{1}(\Omega)}\frac{\|P_{h}Ag-Ag\|_{a}}{\|g\|_{a}}\nonumber\\
&\lesssim & \sup\limits_{g\in
H^{1}(\Omega)}\frac{h^{r}\|Ag\|_{1+r}}{\|g\|_{a}}\nonumber\\
&\lesssim & \sup\limits_{g\in
H^{1}(\Omega)}\frac{h^{r}\|g\|_{a}}{\|g\|_{a}}=h^{r}\rightarrow 0
~(h\rightarrow 0).
\end{eqnarray}
It is clear that $A_{h}$ is a finite rank operator, then $A$ is a completely continuous operator.\\
From \cite{babuska} and \cite{bramble} we know that (\ref{s2.2}) and
(\ref{s2.3}) have the following equivalent operator forms,
respectively:
\begin{eqnarray}\label{s2.11}
&&Au=\mu u,\\\label{s2.12}
&&A_{h}u_{h}=\mu_{h} u_{h},
\end{eqnarray}
where $\mu=\frac{1}{\lambda}, \mu_{h}=\frac{1}{\lambda_{h}}$.
In this paper, $\mu$ and $\mu_{h}$, $\lambda$ and $\lambda_{h}$ are all called eigenvalues.\\
\indent Suppose that $\mu_{k}$ is the $kth$ eigenvalue of $A$ and the algebraic multiplicity of $\mu_{k}$ is equal to
$q$, $\mu_{k}=\mu_{k+1}=\cdots=\mu_{k+q-1}$. Let $M(\mu_{k})$ be the
space spanned by all eigenfunctions corresponding to $\mu_{k}$ of
$A$ and  $M_{h}(\mu_{k})$ be the direct sum of eigenspaces
corresponding to all eigenvalues of $A_{h}$ that converge to
$\mu_{k}$. Let $\widehat{M}(\mu_{k})=\{v: v\in M(\mu_{k}), \|v\|_{a}=1\}$,
$\widehat{M}_{h}(\mu_{k})=\{v: v\in M_{h}(\mu_{k}), \|v\|_{a}=1\}$. We also
write $M(\lambda_{k})=M(\mu_{k})$,
$M_{h}(\lambda_{k})=M_{h}(\mu_{k})$,
$\widehat{M}(\lambda_{k})=\widehat{M}(\mu_{k})$ and
$\widehat{M}_{h}(\lambda_{k})=\widehat{M}_{h}(\mu_{k})$.\\
\indent Denote
\begin{eqnarray*}
&&\sigma(h)=\sup\limits_{f\in H^{1}(\Omega), \|f\|_{a}=1}\inf\limits_{v\in V_{h}}\|Af-v\|_{a},\\
&&\rho(h)=\sup\limits_{f\in L_{2}(\partial\Omega),
\|f\|_{0,\partial\Omega}=1}\inf\limits_{v\in V_{h}}\|Af-v\|_{a},\\
&&\delta_{h}(\lambda_{k})=\sup\limits_{u\in
\widehat{M}(\lambda_{k})}\inf\limits_{v\in V_{h}}\|u-v\|_{a}.
\end{eqnarray*}
It is obvious that
\begin{eqnarray}\label{s2.23}
\delta_{h}(\lambda_{k})\lesssim  \sigma(h)\lesssim  \rho(h).
\end{eqnarray}

The following a priori error estimates have been obtained in \cite{armentano2,liq}:\\
\indent {\bf Lemma 2.2.~} Let $\lambda_{k,h}$ and $\lambda_{k}$ be
the $kth$ eigenvalue of (\ref{s2.3}) and (\ref{s2.2}), respectively.
Then
\begin{eqnarray}\label{s2.13}
\lambda_{k}\leq\lambda_{k,h}\leq \lambda_{k}+ C\delta_{h}^{2}(\lambda_{k});
\end{eqnarray}
for any eigenfunction $u_{k,h}$ corresponding to $\lambda_{k,h}$,
satisfying $\|u_{k,h}\|_{a}=1$, there exists $u_{k}\in
M(\lambda_{k})$ such that
\begin{eqnarray}\label{s2.14}
&&\|u_{k}-u_{k,h}\|_{a}\lesssim
\delta_{h}(\lambda_{k});\\\label{s2.16}
&&\|u_{k}-u_{k,h}\|_{0,\partial\Omega} \lesssim
\rho(h)\delta_{h}(\lambda_{k});\\\label{s2.17}
&&\|u_{k}-u_{k,h}\|_{-\frac{1}{2},\partial\Omega} \lesssim
\sigma(h)\delta_{h}(\lambda_{k});
\end{eqnarray}
for any $u_{k}\in \widehat{M}(\lambda_{k})$, there exists $u_{h}\in
M_{h}(\lambda_{k})$ such that
\begin{eqnarray}\label{s2.15}
\|u_{h}-u_{k}\|_{a}\lesssim \delta_{h}(\lambda_{k}).
\end{eqnarray}

The following lemma states a crucial property (but straightforward) of eigenvalue and eigenfunction approximation.\\
\indent {\bf Lemma 2.3.~} Let $(\lambda,u)$ be an eigenpair of
(\ref{s2.2}), then for any $v \in H^{1}(\Omega), \|v\|_{b}\neq 0$,
the Rayleigh quotient $\frac{a(v,v)}{\|v\|_{b}^{2}}$ satisfies
\begin{eqnarray}\label{s2.18}
\frac{a(v,v)}{\|v\|_{b}^{2}}-\lambda=\frac{\|v-u\|_{a}^{2}}{\|v\|_{b}^{2}}-\lambda
\frac{\|v-u\|_{b}^{2}}{\|v\|_{b}^{2}}.
\end{eqnarray}
\indent {\em Proof.~} See, for instance, Lemma 9.1 in \cite{babuska} for details.~~$\Box$\\

\indent Our analysis is based on the following crucial property of
the shifted-inverse iteration in finite element method (see Lemma 4.1 in \cite{yang3}).\\
\indent {\bf Lemma 2.4.}~~ Let $(\mu_{0}, u_{0})$ be an
approximation for $(\mu_{k}, u_{k})$ where $\mu_{0}$ is not an
eigenvalue of $A_{h}$, and $u_{0}\in V_{h}$ with
$\|u_{0}\|_{a}=1$. Suppose that \\
\indent (C1)~$dist(u_{0}, M_{h}(\mu_{k}))\leq \frac{1}{2}$;\\
\indent (C2)~$|\mu_{0}-\mu_{k}|\leq \frac{\rho}{4}$, $|\mu_{j,h}-\mu_{j}|\leq \frac{\rho}{4}$ for
$j=k-1,k,k+q(j\neq 0)$, where
$\rho=\min\limits_{\mu_{j}\not=\mu_{k}}|\mu_{j}-\mu_{k}|$ is
the separation constant of the eigenvalue $\mu_{k}$;\\
\indent (C3)~ $u'\in V_{h}, u_{k}^{h}\in V_{h}$ satisfy
\begin{eqnarray}\label{s2.19}
(\mu_{0}-A_{h})u'=u_{0},~~~u_{k}^{h}=\frac{u'}{\|u'\|_{a}}.
\end{eqnarray}
Then
\begin{eqnarray}\label{s2.20}
dist(u_{k}^h,M_{h}(\mu_{k}))\leq \frac{4}{\rho}\max\limits_{k\leq
j\leq k+q-1}|\mu_{0}-\mu_{j,h}|dist(u_{0}, M_{h}(\mu_{k})).
\end{eqnarray}

\setcounter{section}{2}\setcounter{equation}{0}

\section{A posteriori error estimate for multi-scale discretizations}

In \cite{by1,bi4,li2,yang3}, a multi-scale discretization scheme based on shifted inverse iteration has been established and its a priori error estimate has been proved. In this section, we will discuss the a posteriori error estimates for this multi-scale discretization scheme, which together with the adaptivity leads to the adaptive algorithm in the next section.

\subsection{Multi-scale discretization scheme and a priori error estimates}
Traditional multigrid methods based on the shifted inverse iteration
has already been used for solving a given discretization scheme
(see, e.g., \cite{hack,shaidurov}) while in this paper we apply the
multi-scale discretization scheme based on the shifted inverse
iteration to solve a differential equation directly.

\indent Let $\{\pi_{h_{i}}\}_{1}^{l}$ be a family of regular triangulations,
and $\{V_{h_{i}}\}_{1}^{l}$ be the conforming
finite element spaces defined on $\{\pi_{h_{i}}\}_{1}^{l}$, and let
$\pi_{H}=\pi_{h_{1}}$, $V_{H}=V_{h_{1}}$,
$\pi_{h}=\pi_{h_{l}}$, $V_{h}=V_{h_{l}}$. Assume that the following condition for meshes holds (see Condition 4.3 in \cite{yang3}).\\
\indent {\bf Condition 3.1.}~~Suppose that $\varepsilon\in (0,1)$ be
a given number, there exists $t_{i}\in (1, 3-\varepsilon]$,
$i=2,3,\cdots$
 such that $\delta_{h_{i}}(\lambda)=\delta_{h_{i-1}}(\lambda)^{t_{i}}$
 and $\delta_{h_{i}}(\lambda)\rightarrow 0~(i\rightarrow \infty)$.\\

\indent Note that Condition 3.1 is closely related the saturation assumption for the
approximation of piecewise polynomials. Recently, it was
proved in \cite{hu1,hu2} that the saturation assumption holds for
the quasi-uniform mesh, which is an essential advance.\\

\noindent{\bf Scheme 1 (multi-scale discretization scheme).}\\
\noindent{\bf Step 1.} Solve (\ref{s2.3}) on $\pi_{H}$: find
$\lambda_{k,H}\in \mathbb{R}, u_{k,H}\in V_{H}$ such that
$\|u_{k,H}\|_{a}=1$ and
\begin{eqnarray*}
a(u_{k,H},v)=\lambda_{k,H} b(u_{k,H},v),~~~\forall v\in V_{H}.
\end{eqnarray*}
\noindent{\bf Step 2.} $u_{k}^{h_{1}}\Leftarrow u_{k,H}$, $\lambda_{k}^{h_{1}}\Leftarrow\lambda_{k,H}$, $i\Leftarrow 2$.\\
\noindent{\bf Step 3.} Solve a linear system on the $\pi_{h_i}$:
find $\tilde{u}\in V_{h_{i}}$ such that
\begin{eqnarray*}
a(\tilde{u}, v)-\lambda_{k}^{h_{i-1}}b(\tilde{u},v)=b(u_{k}^{h_{i-1}},v),~~~\forall
v\in V_{h_{i}}.
\end{eqnarray*}
And set $u_{k}^{h_{i}}=\frac{\tilde{u}}{\|\tilde{u}\|_{a}}$.\\
\noindent{\bf Step 4.} Compute the Rayleigh quotient
$$\lambda_{k}^{h_{i}}=\frac{a(u_{k}^{h_{i}},u_{k}^{h_{i}})}{
b(u_{k}^{h_{i}},u_{k}^{h_{i}})}.$$
 \noindent{\bf Step 5.} If $i=l$,
then output $(\lambda_{k}^{h_{l}}, u_{k}^{h_{l}})$, i.e.,
$(\lambda_{k}^{h}, u_{k}^{h})$, stop. Else, $i\Leftarrow i+1$, and
return to Step 3.\\

\indent $(\lambda_{k}^{h}, u_{k}^{h})$ is used as an approximation for the $kth$ eigenpair, $(\lambda_{k}, u_{k})$, of (\ref{s2.2}).\\

The following a priori error estimates can be obtained by using the proof argument in \cite{by1,yang3}.\\
\indent{\bf Lemma 3.1.} Let $(\lambda^{h_{l}},u^{h_{l}})$ be an
approximate eigenpair obtained by Scheme 1.
 Suppose that the
condition 3.1 holds and $h_{1}$, i.e. $H$ is properly small. Then
there exists $u_{k}\in M(\lambda_{k})$ such that the following error
estimates hold:
\begin{eqnarray}\label{s3.1}
&&\|u_{k}^{h_{l}}-u_{k}\|_{a}\lesssim
\delta_{h_{l}}(\lambda_{k}),\\\label{e3.2}
&&|\lambda_{k}^{h_{l}}-\lambda_{k}|\lesssim
\|u_{k}^{h_{l}}-u_{k}\|_{a}^{2},~~~l\geq 2.
\end{eqnarray}

\subsection{The a posteriori error estimates}
There exists many publications on the a posteriori error estimates (see, e.g., \cite{ainsworth,ainsworth2,babuska3,carstensen2,morin,verfurth}).
Especially, \cite{armentano2} studied the a posteriori error estimate of finite element for the Steklov eigenvalue problem.
Here we discuss the a posteriori error estimate of multi-scale discretization for the
Steklov eigenvalue problem.\\
\indent For any element $T\in\pi_{h}$ with diameter $h_{T}$, let $\mathcal {E}_{T}$ denote the set of edges of $T$, and
$$\mathcal {E}=\bigcup\limits_{T\in\pi_{h}}\mathcal {E}_{T}.$$
We decompose $\mathcal {E}=\mathcal {E}_{\Omega}\cup\mathcal
{E}_{\partial \Omega}$ where $\varepsilon_{\Omega}$ and $\mathcal
{E}_{\partial \Omega}$ refer to interior edges and edges on the
boundary $\partial \Omega$, respectively. For each $\ell\in \mathcal
{E}_{\Omega}$, we choose an arbitrary unit normal vector $n_{\ell}$
and denote the two triangles sharing this
edge by $T_{in}$ and $T_{out}$, where $n_{\ell}$ points outwards $T_{in}$.\\
For $v_{h}\in V_{h}$ we set
\begin{eqnarray*}
[[\frac{\partial v_{h}}{\partial n_{l}}]]_{\ell}=\bigtriangledown
(v_{h}|_{T_{out}})\cdot n_{\ell} -\bigtriangledown
(v_{h}|_{T_{in}})\cdot n_{\ell}.
\end{eqnarray*}

\indent Let $I_{h}^{c}:H^{1}(\Omega)\to V_{h}$ be a
Cl$\acute{e}$ment interpolation operator, then from
\cite{armentano2} or Lemma 1.4 in \cite{verfurth} the following
error estimates for $I_{h}^{c}$ are valid:
\begin{eqnarray}\label{s3.3}
&&\|v-I_{h}^{c}v\|_{0,T}\leq Ch_{T}\|
v\|_{1,\widetilde{T}},\\\label{s3.4}
 &&\|v-I_{h}^{c}v\|_{0,\ell}\leq
C|l|^{\frac{1}{2}}\| v\|_{1,\widetilde{\ell}},
\end{eqnarray}
where $\widetilde{T}$ is the union of all elements sharing a vertex with $T$ and
$\widetilde{\ell}$ is the union of all elements sharing a vertex with $\ell$.\\
Let the eigenvectors $\{u_{j,h_{l}}\}_{k}^{k+q-1}$ be an orthonormal
basis of $M_{h_{l}}(\lambda_{k})$, and denote
\begin{eqnarray}\label{s3.5}
u^{*}=\sum\limits_{j=k}^{k+q-1}a(u_{k}^{h_{l}},
u_{j,h_{l}})u_{j,h_{l}}.
\end{eqnarray}
It follows from Lemma 2.2 that there exist $\{u_{j}^{0}\}_{k}^{k+q-1}\subset M(\lambda_{k})$ making $u_{j,h_{l}}-u_{j}^{0}$ satisfy
(\ref{s2.14}), (\ref{s2.16}) and (\ref{s2.17}). Let
\begin{eqnarray}\label{s3.6}
u=\sum\limits_{j=k}^{k+q-1}a(u_{k}^{h_{l}}, u_{j,h_{l}})u_{j}^{0},
\end{eqnarray}
then $u\in M(\lambda_{k})$ and
\begin{eqnarray}\label{s3.6a}
u-u^{*}=\sum\limits_{j=k}^{k+q-1}a(u_{k}^{h_{l}},u_{j,h_{l}})(u_{j}^{0}-u_{j,h_{l}}).
\end{eqnarray}
Let
$$\widehat{\lambda}_{k,h_{l}}=\frac{1}{q}\sum\limits_{j=k}^{k+q-1}\lambda_{j,h_{l}}.$$

For each $\ell\in \mathcal {E}$ we define the jump residual:
\begin{eqnarray}\label{s3.7}
J_{\ell}(u^{*})&=& \left \{
\begin{array}{ll}\frac{1}{2}[[\frac{\partial u^{*}}{\partial n_{l}}]]_{\ell}~~~ \ell\in \mathcal {E}_{\Omega},\\
\widehat{\lambda}_{k,h_{l}}u^{*}-\frac{\partial u^{*}}{\partial
n_{\ell}}~~~\ell\in \mathcal {E}_{\partial \Omega}.
\end{array}
\right.
\end{eqnarray}
\begin{eqnarray}\label{s3.7a}
J_{\ell}(u_{k}^{h_{l}})&=& \left \{
\begin{array}{ll}\frac{1}{2}[[\frac{\partial u_{k}^{h_{l}}}{\partial n_{\ell}}]]_{\ell}~~~ \ell\in \mathcal {E}_{\Omega},\\
\lambda_{k}^{h_{l}}u_{k}^{h_{l}}-\frac{\partial
u_{k}^{h_{l}}}{\partial n_{\ell}}~~~\ell\in \mathcal {E}_{\partial
\Omega}.
\end{array}
\right.
\end{eqnarray}
Now, the local error indicator is defined as
\begin{eqnarray}\label{s3.8}
\eta_{T}(v)=(h_{T}^{2}\|v\|_{0,T}^{2} +\sum\limits_{\ell\in
\varepsilon_{T}}|\ell|\|J_{\ell}(v)\|_{0,\ell}^{2} )^{1/2},
\end{eqnarray}
then the global error estimator is given by
\begin{eqnarray}\label{s3.9}
\eta_{\Omega}(v)= (\sum\limits_{T\in \pi_{h}}\eta_{T}(v)^{2})^{1/2}.
\end{eqnarray}

We split the error $u-u_{k}^{h_{l}}=(u-u^{*})+(u^{*}-u_{k}^{h_{l}})$. Next, we will estimate the first term $e=u-u^{*}$.\\
\indent We know that for any $v\in V_{h_{l}
}$
\begin{eqnarray}\label{s3.10}
&&\int\limits_{\Omega}\nabla u^*\cdot\nabla v+\int\limits_{\Omega}
u^* v\nonumber\\
&&~~~=\int\limits_{\Omega}\nabla
(\sum\limits_{j=k}^{k+q-1}a(u_{k}^{h_{l}},
u_{j,h_{l}})u_{j,h_{l}})\cdot\nabla v+\int\limits_{\Omega}
\sum\limits_{j=k}^{k+q-1}a(u_{k}^{h_{l}}, u_{j,h_{l}})u_{j,h_{l}} v\nonumber\\
&&~~~=\sum\limits_{j=k}^{k+q-1} a(u_{k}^{h_{l}},
u_{j,h_{l}})(\int\limits_{\Omega}\nabla u_{j,h_{l}}\cdot\nabla
v+\int\limits_{\Omega} u_{j,h_{l}} v)\nonumber\\
&&~~~=\sum\limits_{j=k}^{k+q-1} a(u_{k}^{h_{l}},
u_{j,h_{l}})\lambda_{j,h_{l}}\int\limits_{\partial\Omega}  u_{j,h_{l}} v\nonumber\\
&&~~~=\sum\limits_{j=k}^{k+q-1} a(u_{k}^{h_{l}}, u_{j,h_{l}})(
\lambda_{j,h_{l}}-\widehat{\lambda}_{k,h_{l}})\int\limits_{\partial\Omega}
u_{j,h_{l}} v +\widehat{\lambda}_{k,h_{l}}\int\limits_{\partial\Omega} u^{*}
v,
\end{eqnarray}
thus the error $e=u-u^{*}$ satisfies
\begin{eqnarray}\label{s3.11}
&&\int\limits_{\Omega}\nabla e\cdot\nabla v+\int\limits_{\Omega} e
v=\int\limits_{\partial\Omega}\lambda_{k}u v-\int\limits_{\Omega}\nabla
u^*\cdot\nabla v-\int\limits_{\Omega} u^* v\nonumber\\
&&~~~=\int\limits_{\partial\Omega}\lambda_{k}u
v-\widehat{\lambda}_{k,h_{l}}\int\limits_{\partial\Omega} u^{*} v\nonumber\\
&&~~~~-\sum\limits_{j=k}^{k+q-1} a(u_{k}^{h_{l}}, u_{j,h_{l}})(
\lambda_{j,h_{l}}-\widehat{\lambda}_{k,h_{l}})\int\limits_{\partial\Omega}u_{j,h_{l}} v,~~~\forall v\in V_{h_{l}}.
\end{eqnarray}
On the other hand, for any $v\in H^{1}(\Omega)$, we have
\begin{eqnarray*}
&&\int\limits_{\Omega}\nabla e\cdot\nabla v+\int\limits_{\Omega} e
v=a(u,v)-a(u^*,v)= \int\limits_{\partial\Omega}\lambda_{k} u
v-\sum\limits_{T}\{\int\limits_{\partial T}\frac{\partial
u^*}{\partial n}v+\int\limits_{T}u^{*}v\},
\end{eqnarray*}
and so
\begin{eqnarray}\label{s3.12}
&&\int\limits_{\Omega}\nabla e\cdot\nabla v+\int\limits_{\Omega} e
v\nonumber\\
&&~~~= \sum\limits_{T}\{-\int\limits_{T} u^* v +\sum\limits_{\ell\in
\mathcal {E}_{T}\cap\mathcal {E}_{\partial
\Omega}}\int\limits_{\ell}(\widehat{\lambda}_{k,h_{l}}
u^*-\frac{\partial u^*}{\partial n_{\ell}})v+\frac{1}{2}
\sum\limits_{\ell\in \mathcal {E}_{T}\cap\mathcal {E}_{\Omega}}
\int\limits_{\ell}[[\frac{\partial u^{*}}{\partial
n_{\ell}}]]_{\ell}v\}\nonumber\\
&&~~~~~~ +\int\limits_{\partial\Omega}\lambda_{k} u
v-\int\limits_{\partial\Omega}\widehat{\lambda}_{k,h_{l}} u^*
v,~~\forall v\in H^{1}(\Omega).
\end{eqnarray}

The following lemma provides the global upper bound of $e$.\\
\indent{\bf Lemma 3.2.} The error $e=u-u^{*}$ satisfies
\begin{eqnarray}\label{s3.13}
\|e\|_{a}\lesssim \eta_{\Omega}(u^{*}) +
\sigma(h_{l})\delta_{h_{l}}(\lambda_{k}).
\end{eqnarray}
\indent{\em Proof.}~~By (\ref{s3.11}) and (\ref{s3.12}), noting that the definition of $J_{\ell}(u^{*})$, we deduce that
\begin{eqnarray}\label{s3.14}
&&\int\limits_{\Omega}\nabla e\cdot\nabla e+\int\limits_{\Omega} e
e\nonumber\\
&&~~~= \int\limits_{\Omega}\nabla e\cdot(\nabla e-\nabla I_{h}^{c}e)
+\int\limits_{\Omega} e (e-I_{h}^{c}e)+\int\limits_{\partial\Omega}\lambda_{k}
uI_{h}^{c}e-\int\limits_{\partial\Omega}\widehat{\lambda}_{k,h_{l}}
u^{*}I_{h}^{c}e\nonumber\\
&&~~~~~~- \sum\limits_{j=k}^{k+q-1} a(u_{k}^{h_{l}}, u_{j,h_{l}})(
\lambda_{j,h_{l}}-\widehat{\lambda}_{k,h_{l}})\int\limits_{\partial\Omega}
u_{j,h_{l}} I_{h}^{c}e\nonumber\\
&&~~~= \sum\limits_{T}\{-\int\limits_{T}u^{*}( e-
I_{h}^{c}e)+\sum\limits_{\ell\in
\mathcal {E}_T}\int\limits_{\ell}J_{\ell}(u^{*})
(e-I_{h}^{c}e)\}+\int\limits_{\partial\Omega}(\lambda_{k} u-
\widehat{\lambda}_{k,h_{l}}
u^{*})e\nonumber\\
&&~~~~~~- \sum\limits_{j=k}^{k+q-1} a(u_{k}^{h_{l}}, u_{j,h_{l}})(
\lambda_{j,h_{l}}-\widehat{\lambda}_{k,h_{l}})\int\limits_{\partial\Omega}
u_{j,h_{l}} I_{h}^{c}e.
\end{eqnarray}
Then
\begin{eqnarray}\label{s3.15}
&&\|e\|_{a}^{2} \leq  \sum\limits_{T}\|u^{*}\|_{0,T}\| e-
I_{h}^{c}e\|_{0,T}+ \sum\limits_{T}\sum\limits_{\ell\in \mathcal {E}_T}
\|J_{\ell}(u^{*})\|_{0,\ell}\|e-I_{h}^{c}e\|_{0,\ell}\nonumber\\
&&~~~~~~+|\int\limits_{\partial\Omega}(\lambda_{k} u-
\widehat{\lambda}_{k,h_{l}} u^{*})e|+ |\sum\limits_{j=k}^{k+q-1}
a(u_{k}^{h_{l}}, u_{j,h_{l}})(
\lambda_{j,h_{l}}-\widehat{\lambda}_{k,h_{l}})\int\limits_{\partial\Omega}
u_{j,h_{l}} I_{h}^{c} e|\nonumber\\
&&~~~\lesssim \sum\limits_{T}h_{T}\|u^{*}\|_{0,T}\|
e\|_{1,\widetilde{T}}+ \sum\limits_{T}\sum\limits_{\ell\in
\mathcal {E}_T}|\ell|^{\frac{1}{2}}
\|J_{\ell}(u^{*})\|_{0,\ell}\|e\|_{1,\widetilde{\ell}}
\nonumber\\
&&~~~~~~+|\int\limits_{\partial\Omega}(\lambda_{k} u-
\widehat{\lambda}_{k,h_{l}}
u^{*}) e|+ |\sum\limits_{j=k}^{k+q-1} a(u_{k}^{h_{l}}, u_{j,h_{l}})(
\lambda_{j,h_{l}}-\widehat{\lambda}_{k,h_{l}})\int\limits_{\partial\Omega}
u_{j,h_{l}} I_{h}^{c} e|\nonumber\\
&&~~~\lesssim\{ \sum\limits_{T}(h_{T}^{2}\|u^{*}\|_{0,T}^{2}+
\sum\limits_{\ell\in \mathcal {E}_T}|\ell|
\|J_{\ell}(u^{*})\|_{0,\ell}^{2})\}^{\frac{1}{2}}\|e\|_{a} \nonumber\\
&&~~~~~~+ \|\lambda_{k}u- \widehat{\lambda}_{k,h_{l}} u^{*}\|_{-\frac{1}{2},
\partial\Omega}\|e\|_{a}+\sum\limits_{j=k}^{k+q-1}|\lambda_{j,h_{l}}-\widehat{\lambda}_{k,h_{l}}|\|e\|_{a}.
\end{eqnarray}
By Lemma 2.2, we have
\begin{eqnarray}\label{s3.16}
\|\lambda_{k} u- \widehat{\lambda}_{k,h_{l}} u^{*}\|_{-\frac{1}{2},
\partial\Omega}&\lesssim& |\lambda_{k} - \widehat{\lambda}_{k,h_{l}}|+\|u-
u^{*}\|_{-\frac{1}{2}, \partial\Omega}\nonumber\\
&\lesssim& \sigma(h_{l})\delta_{h_{l}}(\lambda_{k})\\\label{s3.17}
|\lambda_{j,h_{l}}-\widehat{\lambda}_{k,h_{l}}|&\lesssim& \delta_{h_{l}}(\lambda_{k})^{2}.
\end{eqnarray}
Combining (\ref{s3.15}), (\ref{s3.16}) and (\ref{s3.17}), the proof
concludes.~~~$\Box$\\

\indent Next, we shall analyze the local lower bound of $e$.\\
\indent{\bf Lemma 3.3.}~~The error $e=u-u^{*}$ satisfies\\
\indent (a)~ For $T\in\pi_{h}$, if $\partial T\cap\partial
\Omega=\emptyset$, then
\begin{eqnarray}\label{s3.18}
\eta_{T}(u^{*})\lesssim \|e\|_{1,T^{*}},
\end{eqnarray}
\indent where $T^{*}$ denote the union of $T$ and the triangles
sharing an edge with $T$.\\
 \indent (b)~For $T\in \pi_{h}$, if $\partial T\cap\partial \Omega\not=\emptyset$, then
\begin{eqnarray}\label{s3.19}
\eta_{T}(u^{*})\lesssim \|e\|_{1,T}+\sum\limits_{\ell\in \mathcal
{E}_{T}\cap\mathcal {E}_{\partial
\Omega}}|\ell|^{\frac{1}{2}}\|\lambda_{k} u-
\widehat{\lambda}_{k,h_{l}} u^{*}\|_{0,\ell}.
\end{eqnarray}
\indent{\em Proof.}~~First, using the same argument of Lemma 3.1 in  \cite{armentano2} and replacing $\lambda_{h}, u_{h}$ and (3.8) in Lemma 3.1 in \cite{armentano2} by $\widehat{\lambda}_{k,h_{l}}, u^{*}$ and (\ref{s3.12}) in this paper, respectively,
we can prove that
\begin{eqnarray}\label{s3.20}
h_{T}\|u^{*}\|_{0,T}\lesssim \|\nabla e\|_{0,T}+h_{T}\|e\|_{0,T}.
\end{eqnarray}
Also using the same proof method with that of Lemma 3.2 in
\cite{armentano2} and replacing $\lambda_{h}, u_{h}$ and (3.8) in
\cite{armentano2} by $\widehat{\lambda}_{k,h_{l}}, u^{*}$ and
(\ref{s3.12}) here, respectively, we can obtain
\begin{eqnarray}\label{s3.21}
|\ell|^{\frac{1}{2}}\|J_{\ell}(u^{*})\|_{0,\ell}\lesssim(1+h_{T})\|e\|_{1,T}^{2}+|\ell|^{\frac{1}{2}}\|\lambda_{k}
u-\widehat{\lambda}_{k,h_{l}}u^{*}\|_{0,\ell},~~\ell\in\mathcal {E}_{T}
\cap\mathcal {E}_{\partial \Omega},\nonumber\\
\end{eqnarray}
and
\begin{eqnarray}\label{s3.22}
|\ell|^{\frac{1}{2}}\|J_{\ell}(u^{*})\|_{0,\ell}\lesssim
\|e\|_{1,T_{\ell}^{1}\cup T_{\ell}^{2}},~~\ell\in\mathcal {E}_{T}
\cap\mathcal {E}_{\Omega},
\end{eqnarray}
where $T_{\ell}^{1},T_{\ell}^{2}\in \pi_{h}$ are the two triangles
sharing $\ell$.\\
Then, we get (\ref{s3.18}) immediately by combining (\ref{s3.20})
with (\ref{s3.22}), and obtain (\ref{s3.19}) by
 combining (\ref{s3.20}) with (\ref{s3.21}).~~~$\Box$\\

Now, we will analyze the error $u_{k}^{h_{l}}-u^{*}$.\\
\indent{\bf Theorem 3.1.}~~Suppose that the conditions of Lemma
3.1 are satisfied, then
\begin{eqnarray}\label{s3.23}
\|u_{k}^{h_{l}}-u^{*}\|_{a}\lesssim
\delta_{h_{l}}(\lambda_{k})^{\frac{3}{t_{l}}}.
\end{eqnarray}
\indent{\em Proof.} We use Lemma 2.4 to complete the proof. Using
the arguments in \cite{yang3} it is easy to verify that the
conditions
of Lemma 2.4 are satisfied.\\
Using the triangle inequality and (\ref{s2.15}), we get
\begin{eqnarray}\label{s3.26}
dist(u_{0},M_{h_{l}}(\lambda_{k}))\leq dist(u_{0},\widehat{M}(\lambda_{k}))+C\delta_{h_{l}}(\lambda_{k}).
\end{eqnarray}
It follows from (\ref{s2.13}) that $\lambda_{k,h_{l}}\rightarrow
\lambda_{k} (h_{l}\rightarrow 0)$, then by the assumption we have
\begin{eqnarray}\label{s3.27}
&&|\mu_{0}-\mu_{j,h_{l}}|
=|\frac{\lambda_{k}^{h_{l-1}}-\lambda_{j,h_{l}}}{\lambda_{j,h_{l}}\lambda_{k}^{h_{l-1}}}|
=|\frac{\lambda_{k}^{h_{l-1}}-\lambda_{k}+\lambda_{k}-\lambda_{j,h_{l}}}{\lambda_{j,h_{l}}\lambda_{k}^{h_{l-1}}}|\nonumber\\
&&~~~\lesssim
\delta_{h_{l-1}}^{2}(\lambda_{k})+\delta_{h_{l}}^{2}(\lambda_{k}),~~~j=k,k+1,\cdots,
k+q-1.
\end{eqnarray}
\indent Substituting (\ref{s3.26}) and (\ref{s3.27}) into
(\ref{s2.20}), we obtain
\begin{eqnarray}\label{s3.36}
dist(u_{k}^{h_{l}},M_{h_{l}}(\lambda_{k}))\lesssim
(\delta_{h_{l-1}}^{2}(\lambda_{k})+\delta_{h_{l}}^{2}(\lambda_{k}))(dist(u_{0},\widehat{M}(\lambda_{k}))+C\delta_{h_{l}}(\lambda_{k})).
\end{eqnarray}
Observing that
$$dist(u_{k}^{h_{l}},
M_{h_{l}}(\lambda_{k}))=\|u_{k}^{h_{l}}-\sum\limits_{j=k}^{k+q-1}a(u_{k}^{h_{l}},
u_{j,h_{l}})u_{j,h_{l}}\|_{a},$$
and recalling that $u^{*}=\sum\limits_{j=k}^{k+q-1}a(u_{k}^{h_{l}}, u_{j,h_{l}})u_{j,h_{l}}$, then
\begin{eqnarray}\label{s3.37}
\|u_{k}^{h_{l}}-u^{*}\|_{a}& \lesssim&
(\delta_{h_{l-1}}^{2}(\lambda_{k})+\delta_{h_{l}}^{2}(\lambda_{k}))(dist(u_{0},\widehat{M}(\lambda_{k}))+C\delta_{h_{l}}(\lambda_{k}))\nonumber\\
&\lesssim &\delta_{h_{l-1}}(\lambda_{k})^{3} \lesssim
\delta_{h_{l}}(\lambda_{k})^{\frac{3}{t_{l}}}.
\end{eqnarray}
The proof is completed.~~~$\Box$\\

\indent{\bf Lemma 3.4.}~~Suppose that the conditions of Lemma 3.1
are satisfied, then
\begin{eqnarray}\label{s3.38a}
|\eta_{T}(u^{*})-\eta_{T}(u_{k}^{h_{l}})| &\lesssim&
\delta_{h_{l}}(\lambda_{k})^{2}\|u^{*}\|_{1,T}+\|u_{k}^{h_{l}}-u^{*}\|_{1,T},\\\label{s3.38b}
|\eta_{\Omega}(u^{*})-\eta_{\Omega}(u_{k}^{h_{l}})| &\lesssim&
\delta_{h_{l}}(\lambda_{k})^{2}\|u^{*}\|_{a}+\|u_{k}^{h_{l}}-u^{*}\|_{a}.
\end{eqnarray}
\indent{\em Proof.}~~From (\ref{s3.8}) and the triangle inequality, we
deduce
\begin{eqnarray*}
|\eta_{T}(u^{*})-\eta_{T}(u_{k}^{h_{l}})| \leq
(h_{T}^{2}\|u^{*}-u_{k}^{h_{l}}\|_{0,T}^{2} +\sum\limits_{\ell\in
\mathcal {E}_{T}}|\ell|\|J_{\ell}(u^{*})-J_{\ell}(u_{k}^{h_{l}})\|_{0,\ell}^{2}
)^{1/2}.
\end{eqnarray*}
From  (\ref{s3.7}) and  (\ref{s3.7a}), if $\ell\in \mathcal {E}_{\Omega}$, we have
\begin{eqnarray*}
\|J_{\ell}(u^{*})-J_{\ell}(u_{k}^{h_{l}})\|_{0,\ell}=
\frac{1}{2}\|[[\frac{\partial (u^{*}-u_{k}^{h_{l}})}{\partial
n_{\ell}}]]_{\ell}\|_{0,\ell}\lesssim
h_{T}^{-\frac{1}{2}}\|u_{k}^{h_{l}}-u^{*}\|_{1,T},
\end{eqnarray*}
and if $\ell\in \mathcal {E}_{\partial \Omega}$,
\begin{eqnarray*}
&&\|J_{\ell}(u^{*})-J_{\ell}(u_{k}^{h_{l}})\|_{0,\ell}\leq
\|\widehat{\lambda}_{k,h_{l}}u^{*}-\lambda_{k}^{h_{l}}u_{k}^{h_{l}}-\frac{\partial(u^{*}-u_{k}^{h_{l}})}{\partial
n_{\ell}}\|_{0,\ell}\\
&&~~~\lesssim
|\widehat{\lambda}_{k,h_{l}}-\lambda_{k}^{h_{l}}|\|u^{*}\|_{1,T}+\|u_{k}^{h_{l}}-u^{*}\|_{1,T}+h_{T}^{-\frac{1}{2}}\|u_{k}^{h_{l}}-u^{*}\|_{1,T}\\
&&~~~\lesssim
\delta_{h_{l}}(\lambda_{k})^{2}\|u^{*}\|_{1,T}+\|u_{k}^{h_{l}}-u^{*}\|_{1,T}+h_{T}^{-\frac{1}{2}}\|u_{k}^{h_{l}}-u^{*}\|_{1,T}.
\end{eqnarray*}
Combining the above three inequalities we obtain (\ref{s3.38a}). From the triangle inequality and (\ref{s3.38a}) we derive (\ref{s3.38b}).~~~$\Box$\\

Combining Lemmas 3.2-3.4 and Theorem 3.1, we give the global bound and the local lower bound of the error.\\
\indent{\bf Theorem 3.2.}~~Suppose that the conditions of Lemma 3.1
are satisfied, then there exists $u_{k}\in M(\lambda_{k})$ such that
\begin{eqnarray}\label{s3.38}
\|u_{k}-u_{k}^{h_{l}}\|_{a} \lesssim \eta_{\Omega}(u_{k}^{h_{l}}).
\end{eqnarray}
\indent{\em Proof.}~~Select $u_{k}\in M(\lambda_{k})$ which is given
by (\ref{s3.6}), then from Lemma 3.2, Theorem 3.1 and Lemma 3.4 we
get
\begin{eqnarray*}
\|u_{k}-u_{k}^{h_{l}}\|_{a}& \leq&
\|u_{k}-u^{*}\|_{a}+\|u^{*}-u_{k}^{h_{l}}\|_{a}\nonumber\\
 &\lesssim& \eta_{\Omega}(u^{*}) +
\sigma(h_{l})\delta_{h_{l}}(\lambda_{k})+
\delta_{h_{l}}(\lambda_{k})^{\frac{3}{t_{l}}}\\
&\lesssim& \eta_{\Omega}(u_{k}^{h_{l}}) +
\sigma(h_{l})\delta_{h_{l}}(\lambda_{k})+
\delta_{h_{l}}(\lambda_{k})^{\frac{3}{t_{l}}}.
\end{eqnarray*}
Since in the above inequality the first term on the right hand side is the dominant term and the other two terms are of high order,
(\ref{s3.38}) is true.~~~$\Box$\\

\indent{\bf Theorem 3.3.}~~Suppose that the conditions of Lemma 3.1
are satisfied, then there exists $u_{k}\in M(\lambda_{k})$ such that
\begin{eqnarray}\label{s3.39}
\eta_{\Omega}(u_{k}^{h_{l}}) \lesssim \|u_{k}-u_{k}^{h_{l}}\|_{a}.
\end{eqnarray}
\indent{\em Proof.}~~Select $u_{k}\in M(\lambda_{k})$ which is given
by (\ref{s3.6}), then from Lemma 3.3, (\ref{s3.6a}),  (\ref{s2.13})
and  (\ref{s2.16}) we have
\begin{eqnarray*}
&&\sum\limits_{T\in\pi_{h_{l}}}\eta_{T}(u^{*})^{2} \lesssim
\sum\limits_{T\in\pi_{h_{l}}}\|e\|_{1,T}^{2} + \sum\limits_{\ell\in
\mathcal {E}_{\partial \Omega}}|\ell|\|\lambda_{k} u_{k}-
\widehat{\lambda}_{k,h_{l}}
u^{*}\|_{0,\ell}^{2}\nonumber\\
&&~~~\lesssim \|e\|_{a}^{2} + h_{l}\|\lambda_{k} u_{k}-
\widehat{\lambda}_{k,h_{l}}
u^{*}\|_{0,\partial\Omega}^{2}\nonumber\\
&&~~~\lesssim \|e\|_{a}^{2} + h_{l}(|\lambda_{k} -
\widehat{\lambda}_{k,h_{l}}|^{2}+\|u_{k}-
u^{*}\|_{0,\partial\Omega}^{2})\nonumber\\
&&~~~\lesssim \|e\|_{a}^{2} + h_{l}\rho(h_{l})^{2}\delta_{h_{l}}(\lambda_{k})^{2},
\end{eqnarray*}
thus
\begin{eqnarray*}
\eta_{\Omega}(u^{*})& \lesssim& \|e\|_{a} +
h_{l}^{\frac{1}{2}}\rho(h_{l})\delta_{h_{l}}(\lambda_{k})\\
&\leq&  \|u_{k}-u_{k}^{h_{l}}\|_{a} +\|u_{k}^{h_{l}}-u^{*}\|_{a}+
h_{l}^{\frac{1}{2}}\rho(h_{l})\delta_{h_{l}}(\lambda_{k}),
\end{eqnarray*}
which together with Theorem 3.1 and  Lemma 3.4 leads to
\begin{eqnarray}\label{s3.400}
\eta_{\Omega}(u_{k}^{h_{l}}) \lesssim
\|u_{k}-u_{k}^{h_{l}}\|_{a}+\delta_{h_{l}}(\lambda_{k})^{2}
+\delta_{h_{l}}(\lambda_{k})^{\frac{3}{t_{l}}}
+h_{l}^{\frac{1}{2}}\rho(h_{l})\delta_{h_{l}}(\lambda_{k}).
\end{eqnarray}
Noting that $\frac{3}{t_{l}}>1$, we know the first term on the right hand side of (\ref{s3.400}) is the dominant term and the others are of high order,
then we derive (\ref{s3.39}).~~~$\Box$\\

\indent{\bf Theorem 3.4.}~~Under the conditions of Lemma 3.1, there exists $u_{k}\in M(\lambda_{k})$ such that\\
\indent (a)~ For $T\in\pi_{h_{l}}$, if $\partial T\cap\partial
\Omega=\emptyset$, then
\begin{eqnarray}\label{s3.40}
\eta_{T}(u_{k}^{h_{l}})\lesssim
\|u_{k}-u_{k}^{h_{l}}\|_{1,T^{*}}+\|u_{k}^{h_{l}}-u^{*}\|_{1,T^{*}}+\delta_{h_{l}}(\lambda_{k})^{2}\|u^{*}\|_{1,T},
\end{eqnarray}
\indent where $T^{*}$ denote the union of $T$ and the triangles
sharing an edge with $T$.\\
 \indent (b)~For $T\in \pi_{h_{l}}$, if $\partial T\cap\partial \Omega\not=\emptyset$, then
\begin{eqnarray}\label{s3.41}
&&\eta_{T}(u_{k}^{h_{l}})\lesssim
\|u_{k}-u_{k}^{h_{l}}\|_{1,T}+\|u_{k}^{h_{l}}-u^{*}\|_{1,T}\nonumber\\
&&~~~~~~+\delta_{h_{l}}(\lambda_{k})^{2}\|u^{*}\|_{1,T}+\sum\limits_{\ell\in
\mathcal {E}_{T}\cap\mathcal {E}_{\partial
\Omega}}|\ell|^{\frac{1}{2}}\|\lambda_{k} u_{k}-
\widehat{\lambda}_{k,h_{l}} u^{*}\|_{0,\ell}.
\end{eqnarray}
\indent{\bf Proof.}~~Select $u_{k}\in M(\lambda_{k})$ which is given by (\ref{s3.6}), then from the triangle inequality we have
$$\|u_{k}-u^{*}\|_{1,T^{*}}\leq \|u_{k}-u_{k}^{h_{l}}\|_{1,T^{*}}+\|u_{k}^{h_{l}}-u^{*}\|_{1,T^{*}}.$$
From (\ref{s3.18}) we obtain
\begin{eqnarray*}
\eta_{T}(u^{*})\lesssim
\|u_{k}-u_{k}^{h_{l}}\|_{1,T^{*}}+\|u_{k}^{h_{l}}-u^{*}\|_{1,T^{*}},
\end{eqnarray*}
combining with Lemma 3.4 we get (\ref{s3.40}).\\
Similarly, from (\ref{s3.19}) and Lemma 3.4 we know that (\ref{s3.41}) is valid.~~~$\Box$\\

\indent{\bf Remark 3.2.}~~Since
$$\|\lambda_{k}
u_{k}- \widehat{\lambda}_{k,h_{l}}
u^{*}\|_{0,\partial\Omega}\lesssim
\lambda_{k}\|u_{k}-u^{*}\|_{0,\partial\Omega}+
|\lambda_{k}-\widehat{\lambda}_{k,h_{l}}|\lesssim
\rho(h_{l})\delta_{h_{l}}(\lambda_{k})+\delta_{h_{l}}^{2}(\lambda_{k}),$$
according to Remark 3.1 in \cite{armentano2} we know that the term
$\sum\limits_{\ell\in \mathcal {E}_{T}\cap\mathcal {E}_{\partial
\Omega}}|\ell|^{\frac{1}{2}}\|\lambda_{k} u_{k}-
\widehat{\lambda}_{k,h_{l}} u^{*}\|_{0,\ell}$ is a higher order
term. Similarly, $\|u_{k}^{h_{l}}-u^{*}\|_{1,T^{*}}$ and
$\|u_{k}^{h_{l}}-u^{*}\|_{1,T}$ are higher order terms. In fact,
from Theorem 3.1 we have
\begin{eqnarray*}
\|u_{k}^{h_{l}}-u^{*}\|_{a}\lesssim
\delta_{h_{l}}(\lambda_{k})^{\frac{3}{t_{l}}}~~~(\frac{3}{t_{l}}>1).
\end{eqnarray*}
And it is obviously that $\delta_{h_{l}}(\lambda_{k})^{2}\|u^{*}\|_{1,T}$ is a higher order term. Therefore, in (\ref{s3.40}) and (\ref{s3.41}) the first
term on the right-hand side is the dominant term.\\

\indent Using the proof method of (4.21) in \cite{yang3} we can
prove that $\|u_{k}^{h_{l}}-u_{k}\|_{b}$ is an infinitesimal of
higher order comparing with $\|u_{k}^{h_{l}}-u_{k}\|_{a}$, thus,
from Lemma 2.3 it is easy to prove that
$\lambda_{k}^{h_{l}}-\lambda_{k}=\mathcal
{O}(\|u_{k}^{h_{l}}-u_{k}\|_{a}^{2})$. Combining
(\ref{s3.38}) with (\ref{s3.39}) we can obtain the following estimates for approximate eigenvalue.\\
\indent{\bf Theorem 3.5.}~~Suppose that the conditions of Lemma 3.1
are satisfied, then
 \begin{eqnarray}\label{s3.411}
\eta_{\Omega}(u_{k}^{h_{l}})^{2}\lesssim|\lambda_{k}^{h_{l}}-\lambda_{k}| \lesssim \eta_{\Omega}(u_{k}^{h_{l}})^{2}.
\end{eqnarray}

\indent Theorem 3.5 shows that $\eta_{\Omega}(u_{k}^{h_{l}})^{2}$ is a reliable and effective estimator for $\lambda_{k}^{h_{l}}$.

\setcounter{section}{3}\setcounter{equation}{0}

\section{Adaptive algorithms and numerical experiments}
\indent In this section, based on the a posteriori error estimates we will establish an adaptive procedure of shifted inverse iteration type and present another two adaptive methods for the Steklov eigenvalue problems, and report several numerical experiments to compare the efficiency of three different adaptive algorithms and illustrate the efficiency of our adaptive method.\\

\subsection{Adaptive algorithms based on multi-scale discretization}

\indent The following adaptive Algorithm 4.1 is usual and standard, which was discussed in \cite{gm}.\\
\indent{\bf Algorithm 4.1}~~Choose parameter $0<\omega<1$.\\
\indent{\bf Step 1.}~Pick any initial mesh $\pi_{h_{1}}$ with mesh size $h_{1}$.\\
\indent{\bf Step 2.}~Solve (\ref{s2.3}) on $\pi_{h_{1}}$ for discrete solution $(\lambda_{k,h_{1}}, u_{k,h_{1}})$.\\
\indent{\bf Step 3.}~Let $l=1$.\\
\indent{\bf Step 4.}~Compute the local indicators ${\eta}_{T}(u_{k,h_{l}})$.\\
\indent{\bf Step 5.}~~Construct
$\widehat{\pi}_{h_{l}}\subset\pi_{h_{l}}$ by {\bf Marking Strategy
E1}
and parameter $\omega$.\\
\indent{\bf Step 6.}~Refine $\pi_{h_{l}}$ to get a new mesh
$\pi_{h_{l+1}}$
by Procedure ${\bf REFINE}$.\\
\indent{\bf Step 7.}~Solve (\ref{s2.3}) on $\pi_{h_{l+1}}$ for discrete solution $(\lambda_{k,h_{l+1}}, u_{k,h_{l+1}})$.\\
\indent{\bf Step 8.}~Let $l=l+1$ and go to Step 4.\\

\indent {\bf Marking Strategy E1}\\
\indent Given parameter $0<\omega<1$:\\
\indent{\bf Step 1.}~~Construct a minimal subset
$\widehat{\pi}_{h_{l}}$ of $\pi_{h_{l}}$ by selecting some elements
in $\pi_{h_{l}}$ such that
\begin{eqnarray*}
\sum\limits_{T\in \widehat{\pi}_{h_{l}}}{\eta}_{T}(u_{k,h_{l}})^{2}
\geq \omega{\eta}_{\Omega}(u_{k,h_{l}})^{2}.
\end{eqnarray*}
\indent{\bf Step 2.}~~Mark all the elements in
$\widehat{\pi}_{h_{l}}$.\\

In Algorithm 4.1 and   {\bf Marking Strategy E1}, the a posteriori error estimators
${\eta}_{T}(u_{k,h_{l}})$ and ${\eta}_{\Omega}(u_{k,h_{l}})$
are defined by (\ref{s3.8}), (\ref{s3.9}) and (\ref{s3.7a}) with $u_{k}^{h_{l}}$
and $\lambda_{k}^{h_{l}}$ replaced by $u_{k,h_{l}}$ and $\lambda_{k,h_{l}}$, respectively.\\
Based on the work of \cite{dr,mm,rs,pssg}, we give the following adaptive algorithm 4.2 for the Steklov eigenvalue problem.\\
\indent {\bf Algorithm 4.2}~~ Choose parameter $0<\omega<1$.\\
\indent{\bf Step 1.}~~Pick any initial mesh $\pi_{h_{1}}$.\\
\indent{\bf Step 2.}~~Solve (\ref{s2.3}) on $\pi_{h_{1}}$
for discrete solution $(\lambda_{k,h_{1}}, u_{k,h_{1}})$. $u_{k*}^{h_{1}}\Leftarrow u_{k,h_{1}}$, $\lambda_{k*}^{h_{1}}\Leftarrow \lambda_{k,h_{1}}$.\\
\indent{\bf Step 3.}~~Let $l=1$.\\
\indent{\bf Step 4.}~~Compute the local indicators ${\eta}_{T}(u_{k*}^{h_{l}})$.\\
\indent{\bf Step 5.}~~Construct
$\widehat{\pi}_{h_{l}}\subset\pi_{h_{l}}$ by {\bf Marking Strategy
E2}
and parameter $\omega$.\\
\indent{\bf Step 6.}~~Refine $\pi_{h_{l}}$ to get a new mesh
$\pi_{h_{l+1}}$
by Procedure ${\bf REFINE}$.\\
\indent{\bf Step 7.}~~Find $\tilde{u}\in V_{h_{l+1}}$ such that
\begin{eqnarray*}
a(\tilde{u},\psi)=b(u_{k*}^{h_{l}},\psi), ~~~
\forall \psi \in V_{h_{l+1}};
\end{eqnarray*}
Denote $u_{k*}^{h_{l+1}}=\frac{\tilde{u}}{\|\tilde{u}\|_{a}}$ and
compute the Rayleigh quotient
\begin{eqnarray*}
\lambda_{k*}^{h_{l+1}}=\frac{a(u_{k*}^{h_{l+1}},u_{k*}^{h_{l+1}})}{b(u_{k*}^{h_{l+1}},u_{k*}^{h_{l+1}})}.
\end{eqnarray*}
\indent{\bf Step 8.}~Let $l=l+1$ and go to Step 4.\\

\indent {\bf Marking Strategy E2} is the same as {\bf Marking Strategy E1} except that the corresponding a posteriori error estimators
are taken as $\eta_{T}(u_{k*}^{h_{l}})$ and $\eta_{\Omega}(u_{k*}^{h_{l}})$, and $\eta_{T}(u_{k*}^{h_{l}})$ and $\eta_{\Omega}(u_{k*}^{h_{l}})$
are defined by (\ref{s3.8}), (\ref{s3.9}) and (\ref{s3.7a}) with $u_{k}^{h_{l}}$ and $\lambda_{k}^{h_{l}}$ replaced by $u_{k*}^{h_{l}}$
and $\lambda_{k*}^{h_{l}}$, respectively.\\

Combining the a posteriori error estimator and Scheme 1, we establish the following adaptive algorithm.\\
\indent {\bf Algorithm 4.3}~~ Choose parameter $0<\omega<1$.\\
\indent{\bf Step 1.}~~Pick any initial mesh $\pi_{h_{1}}$.\\
\indent{\bf Step 2.}~~Solve (\ref{s2.3}) on $\pi_{h_{1}}$
for discrete solution $(\lambda_{k,h_{1}}, u_{k,h_{1}})$. $u_{k}^{h_{1}}\Leftarrow u_{k,h_{1}}$, $\lambda_{k}^{h_{1}}\Leftarrow \lambda_{k,h_{1}}$.\\
\indent{\bf Step 3.}~~Let $l=1$.\\
\indent{\bf Step 4.}~~Compute the local indicators ${\eta}_{T}(u_{k}^{h_{l}})$.\\
\indent{\bf Step 5.}~~Construct
$\widehat{\pi}_{h_{l}}\subset\pi_{h_{l}}$ by {\bf Marking Strategy
E3}
and parameter $\omega$.\\
\indent{\bf Step 6.}~~Refine $\pi_{h_{l}}$ to get a new mesh
$\pi_{h_{l+1}}$
by Procedure ${\bf REFINE}$.\\
\indent{\bf Step 7.}~~Find $\tilde{u}\in V_{h_{l+1}}$ such that
\begin{eqnarray*}
a(\tilde{u},\psi)-
\lambda_{k}^{h_{l}}b(\tilde{u},\psi)=b(u_{k}^{h_{l}},\psi), ~~~
\forall \psi \in V_{h_{l+1}};
\end{eqnarray*}
Denote $u_{k}^{h_{l+1}}=\frac{\tilde{u}}{\|\tilde{u}\|_{a}}$ and
compute the Rayleigh quotient
\begin{eqnarray*}
\lambda_{k}^{h_{l+1}}=\frac{a(u_{k}^{h_{l+1}},u_{k}^{h_{l+1}})}{b(u_{k}^{h_{l+1}},u_{k}^{h_{l+1}})}.
\end{eqnarray*}
\indent{\bf Step 8.}~Let $l=l+1$ and go to Step 4.\\

\indent  {\bf Marking Strategy E3} is the same as {\bf Marking Strategy E1} except that the corresponding a posteriori error estimators are taken as $\eta_{T}(u_{k}^{h_{l}})$ and $\eta_{\Omega}(u_{k}^{h_{l}})$, and $\eta_{T}(u_{k}^{h_{l}})$ and $\eta_{\Omega}(u_{k}^{h_{l}})$
are defined by (\ref{s3.7a}), (\ref{s3.8}) and (\ref{s3.9}).\\

\subsection{Numerical experiments}

\indent We will report two numerical examples to show the
performances of Algorithms 4.1-4.3. We use MATLAB 2011b to solve
Example 4.1 and Example 4.2. Our program makes use of the
package of Chen \cite{cl}.
We refer to the adaptive program of Chen and take $\omega=0.25$.\\

\indent In step 7 of Algorithm 4.1, we use the internal command eigs in MATLAB to perform the eigenvalue computations. The calling format of eigs we used is eigs(K,M,1,sigma), whose function is to solve an eigenvalue which is closest to sigma. In our computation, K and M are stiffness and mass matrices, respectively, sigma is an eigenvalue derived from the last iteration.\\
\indent For convenience of reading, we specify the following notations used in our tables and figures:\\
\indent ${\it \lambda_{k,h_{l}}}$: The $kth$ eigenvalue derived from the $lth$ iteration obtained by Algorithm 4.1.\\
\indent ${\it N_{k,l}}$: The degrees of freedom of the $lth$ iteration for ${\it \lambda_{k,h_{l}}}$.\\
\indent ${\it CPU_{k,l}(s)}$: The CPU time(s) from the program
starting to the calculate result of the $lth$ iteration
 appearing by using Algorithm 4.1.\\
\indent $\lambda_{k*}^{h_{l}}$: The $kth$ eigenvalue derived from the $lth$ iteration obtained by Algorithm 4.2.\\
\indent $N_{\it k*,l}$: The degrees of freedom of the $lth$ iteration for $\lambda_{k*}^{h_{l}}$.\\
\indent $CPU_{k*,l}(s)$: The CPU time(s) from the program starting to the calculate result of the $lth$ iteration appearing by using Algorithm 4.2.\\
\indent $\lambda_{k}^{h_{l}}$: The $kth$ eigenvalue derived from the $lth$ iteration obtained by Algorithm 4.3.\\
\indent $N_{\it k+,l}$: The degrees of freedom of the $lth$ iteration for $\lambda_{k}^{h_{l}}$.\\
\indent $CPU_{k+,l}(s)$: The CPU time(s) from the program starting to the calculate result of the $lth$ iteration appearing by using Algorithm 4.3.\\
\indent $e_{i}(i=1,2,3)$: The absolute value of the error of approximate eigenvalue obtained by Algorithm $4.i$ $(i=1,2,3)$.\\
\indent $\eta_{i}(i=1,2,3)$: The a posteriori error estimator of approximate eigenvalue obtained by Algorithm $4.i$ $(i=1,2,3)$.\\

\noindent {\bf Example 4.1.} We compute the approximations of the first, the second and the fourth eigenvalue of (\ref{s2.1})
 with the triangle linear finite element on $\bar{\Omega}=[0,1]\times[0,1]$ by Algorithms 4.1,
4.2 and 4.3, respectively, and the results are shown in Tables 1 and 2.
We also depict the error curves
and the a posteriori error estimators of Algorithms 4.1, 4.2 and 4.3 in Figs. 1-3.

\begin{table}[htbp]
\centering \caption{\small The results of Example 4.1 obtained by Algorithms
4.1 and 4.3}
\begin{tabular}{cccccccc}
  \hline
{\it k} & {\it l} & $N_{\it k+,l}$ & $\lambda_{k}^{h_{l}}$ & $CPU_{k+,l}(s)$ & ${\it N_{k,l}}$ & ${\it \lambda_{k,h_{l}}}$ & ${\it CPU_{k,l}(s)}$\\
\hline
1 & 23 & 386414 & 0.24007910 & 237.41 & 386414 & 0.24007910 & 411.07\\
1 & 28 & 776862 & 0.24007909 & 510.34 & ------ & ------------ & -------\\
2 & 27 & 500464 & 1.49230435 & 340.08 & 408337 & 1.49230468 & 553.05\\
2 & 30 & 756961 & 1.49230397 & 531.99 & ------ & ------------ & -------\\
4 & 24 & 404451 & 2.08265532 & 250.73 & 404451 & 2.08265532 & 445.37\\
4 & 28 & 776445 & 2.08265094 & 476.95 & ------ & ------------ & -------\\
  \hline
\end{tabular}
\begin{center}
\footnotesize $\star$ The symbol '------' means that the calculation
by Algorithm 4.1 cannot proceed since the computer runs out of
memory.
\end{center}
%
\centering \caption{\small The results of Example 4.1 obtained by Algorithms
4.2 and 4.3}
\begin{tabular}{cccccccc}
  \hline
{\it k} & {\it l} & $N_{\it k+,l}$ & $\lambda_{k}^{h_{l}}$ & $CPU_{k+,l}(s)$ &  ${\it N_{k*,l}}$ & ${\it \lambda_{k*}^{h_{l}}}$ & ${\it CPU_{k*,l}(s)}$\\
\hline
1 & 23 & 386414 & 0.24007910 & 237.41 &  386267 & 0.24007911 & 199.56 \\
1 & 28 & 776862 & 0.24007909 & 510.34 &  776224 & 0.24007909 & 421.38 \\
2 & 27 & 500464 & 1.49230435 & 340.08 &  261640 & 0.24007911 & 192.71  \\
2 & 30 & 756961 & 1.49230397 & 531.99 &  433359 & 0.24007910 & 276.27  \\
4 & 24 & 404451 & 2.08265532 & 250.73 &  242547 & 0.24007911 & 133.72  \\
4 & 28 & 776445 & 2.08265094 & 476.95 &  434294 & 0.24007910 & 235.98  \\
  \hline
\end{tabular}
\centering \caption{\small The results of Example 4.1 obtained by Algorithm 4.2}
\begin{tabular}{ccccccccc}
  \hline
$l$ & $N_{\it 1*,l}$ & $\lambda_{1*}^{h_{l}}$ & $l$ & $N_{\it 2*,l}$ & $\lambda_{2*}^{h_{l}}$ & $l$ & $N_{\it 4*,l}$ & $\lambda_{4*}^{h_{l}}$\\
\hline
1 & 16641 & 0.24007967 & 1 & 16641 & 1.49234096 & 1 & 16641 & 2.08289558\\
7 & 33775 & 0.24007927 & 17 & 141679 & 1.48357414 & 7 & 32621 & 1.60115184\\
9 & 44671 & 0.24007922 & 18 & 146349 & 1.22503676 & 8 & 32798 & 0.30677585\\
11 & 58295 & 0.24007919& 19 & 146609 & 0.34911707 & 9 & 32945 & 0.24099836\\
16 & 122448 & 0.24007913 & 20 & 146950 & 0.24316279 & 10 & 33637 & 0.24009151\\
28 & 776224 & 0.24007909 & 30 & 433359 & 0.24007910 & 28 & 434294 & 0.24007910\\
  \hline
\end{tabular}
\end{table}

\begin{figure}[htp]
  \centering
  \includegraphics[width=6cm]{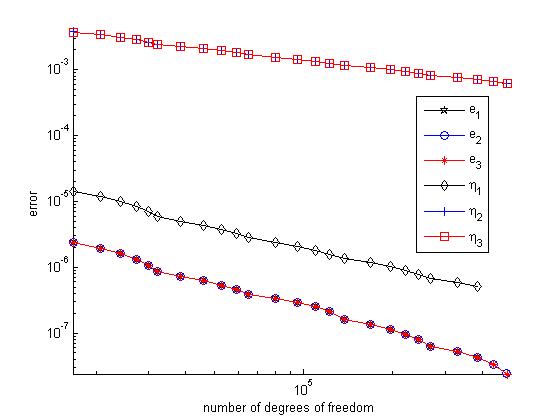}\\
  \caption{\small  The curves of error and the a posteriori estimators for the first eigenvalue of Example 4.1 in log-log scale}
\end{figure}

\begin{figure}[htp]
  \centering
  \includegraphics[width=6cm]{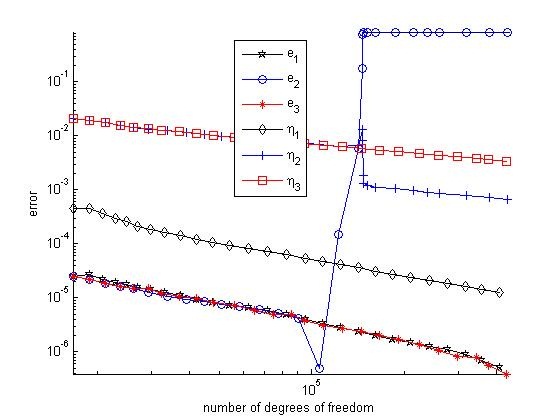}\\
  \caption{\small The curves of error and the a posteriori estimators for the second eigenvalue of Example 4.1 in log-log scale}
\end{figure}

\begin{figure}[htp]
   \centering
  \includegraphics[width=6cm]{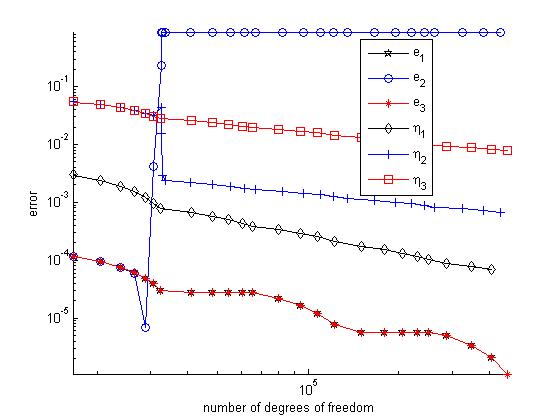}\\
  \caption{\small The curves of error and the a posteriori estimators for the fourth eigenvalue of Example 4.1 in log-log scale}
\end{figure}

\indent We can see from Table 1 that compared with Algorithm 4.1, Algorithm
4.3 costs less CPU time to obtain the same accurate
approximations. We think the reason is although Algorithm 4.1 and 4.3 all use the shifted inverse iteration method,
Step 7 in Algorithm 4.1 is to solve an eigenvalue problem using the command eigs(K,M,1,$\lambda_{k,h_{l}}$), and
the iteration times is usually taken as 2 or 3; while Step 7 in Algorithm 4.3 is to solve a linear system
and the initial iteration value is taken as $u_{k}^{h_{l}}$ which is derived from the last iteration,
and the theoretical analysis assures that there is only one iteration needed.

\indent It can be seen from Table 2 that we can get the same accurate
approximation for the first eigenvalue by Algorithms 4.2 and 4.3,
and Algorithm 4.2 spends less time. However, it is worthy noting in
Table 2 that the approximate eigenvalues $\lambda_{2*}^{h_{l}}$ and
$\lambda_{4*}^{h_{l}}$ obtained by Algorithm 4.2 are not good
approximations for the second and the fourth exact eigenvalue but
are close to the first exact one. We also list some results in Table 3
to show this trend  (In Figs 2 and 3, the curves of $e_2$ also reflect this phenomenon). In fact, from the point of
numerical algebra, we can view Algorithm 4.2 as the inverse
iteration method. And it is well-known that the inverse iteration
method is only applicable to solving the smallest eigenvalue.
\cite{pssg} has pointed out that after an orthogonalization
procedure, Algorithm 4.2 can be suitable to approximate any eigenpair and not only the first one.


In Fig. 1 we can see that the curves of $e_{1}, e_{2}$ and $e_{3}$
are overlapping, and the curves of $\eta_{1}, \eta_{2}$ and
$\eta_{3}$ are basically parallel to those of $e_{1}, e_{2}$ and
$e_{3}$, respectively, which shows that our estimator is reliable
and efficient. In Fig. 2 and 3 we can see that when the numbers of
degrees of freedom are not very large, the curves of $e_{1}, e_{2}$
and $e_{3}$ are overlapping
and the curves of $e_{3}$ and $\eta_{3}$ are basically parallel which illustrate that our estimator is reliable and efficient.\\

\noindent {\bf Example 4.2.} We compute the three smallest
approximate eigenvalues of (\ref{s2.1}) 
with the triangle linear
finite element on $([0,1]\times [0,\frac{1}{2}])\bigcup
([0,\frac{1}{2}]\times [\frac{1}{2},1])$ by Algorithms 4.1, 4.2 and
4.3, respectively, and the results are listed in Table 4 and 5. We
generate the initial mesh by a uniform triangulation with the
diameter $\frac{\sqrt{2}}{128}$, and in Figs. 4-6 we shows the
adaptively refined meshes for $\lambda_{i}(i=1,2,3)$.
We also depict the error curves and the a posteriori error estimators of Algorithms 4.1, 4.2 and 4.3 in Figs. 7-9.

\begin{table}[htbp]
\centering \caption{\small The results for the three smallest
eigenvalues of Example 4.2 obtained by Algorithms 4.1 and 4.3}
\begin{tabular}{cccccccc}
  \hline
{\it k} & {\it l} & $N_{\it k+,l}$ & $\lambda_{k}^{h_{l}}$ & $CPU_{k+,l}(s)$ & ${\it N_{k,l}}$ & ${\it \lambda_{k,h_{l}}}$ & ${\it CPU_{k,l}(s)}$\\
\hline
1 & 25 & 354231 & 0.18296426 & 213.98 & 354231 & 0.18296426 & 383.75\\
1 & 31 & 785693 & 0.18296424 & 533.99 & ------ & ------------ & -------\\
2 & 38 & 433695 & 0.89364798 & 285.66 & 433695 & 0.89364798 & 496.74\\
2 & 42 & 768861 & 0.89364690 & 525.16 & ------ & ------------ & -------\\
3 & 27 & 406124 & 1.68860273 & 237.45 & 406124 & 1.68860273 & 441.66\\
3 & 32 & 801368 & 1.68860181 & 511.28 & ------ & ------------ & -------\\
  \hline
\end{tabular}
%
\centering \caption{\small The results for the three smallest
eigenvalues of Example 4.2 obtained by Algorithms 4.2 and 4.3}
\begin{tabular}{cccccccc}
  \hline
{\it k} & {\it l} & $N_{\it k+,l}$ & $\lambda_{k}^{h_{l}}$ & $CPU_{k+,l}(s)$ & ${\it N_{k*,l}}$ & ${\it \lambda_{k*}^{h_{l}}}$ & ${\it CPU_{k*,l}(s)}$\\
\hline
1 & 25 & 354231 & 0.18296426 & 213.98 &  354231 & 0.18296426 & 192.43  \\
1 & 31 & 785693 & 0.18296424 & 533.99 &  785693 & 0.18296424 & 460.75  \\
2 & 38 & 433695 & 0.89364798 & 285.66 &  394017 & 0.18296425 & 225.75  \\
2 & 42 & 768861 & 0.89364690 & 525.16 &  716738 & 0.18296425 & 402.61  \\
3 & 27 & 465753 & 1.68860256 & 237.45 &  268505 & 0.18296426 & 118.36  \\
3 & 32 & 801368 & 1.68860181 & 511.28 &  445339 & 0.18296425 & 246.42  \\
  \hline
\end{tabular}
%
\centering \caption{\small The results for the three smallest
eigenvalues of Example 4.2 obtained by Algorithm 4.2}
\begin{tabular}{ccccccccc}
  \hline
$l$ & $N_{\it 1*,l}$ & $\lambda_{1*}^{h_{l}}$ & $l$ & $N_{\it 2*,l}$ & $\lambda_{2*}^{h_{l}}$ & $l$ & $N_{\it 3*,l}$ & $\lambda_{3*}^{h_{l}}$\\
\hline
1 & 12545 & 0.18296492 & 1 & 12545 & 0.89423511 & 1 & 12545 & 1.68870013\\
6 & 24006 & 0.18296450 & 9 & 12719 & 0.88859988 & 6 & 20431 & 1.68485494\\
12 & 56774 & 0.18296434 & 11 & 13133 & 0.32003522 & 8 & 22958 & 0.26100111\\
18 & 132597 & 0.18296429 & 13 & 13519 & 0.18326301 & 10 & 24108 & 0.18297591\\
25 & 354231 & 0.18296426 & 38 & 394017 & 0.18296425& 27 & 225747 & 0.18296426\\
31 & 785693 & 0.18296424 & 42 & 716738 & 0.18296425& 32 & 445339 & 0.18296425\\
  \hline
\end{tabular}
\end{table}

\begin{figure}[htp]
  \centering
  \includegraphics[width=4cm]{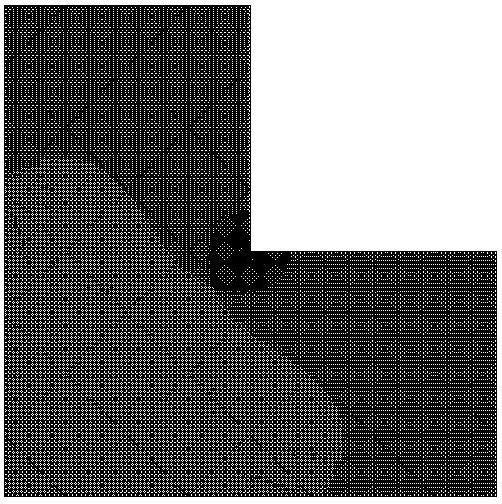}\\
  \caption{\small The adaptively refined mesh with 27113 degrees of freedom for $\lambda_1$ on $([0,1]\times [0,\frac{1}{2}])\bigcup ([0,\frac{1}{2}]\times
[\frac{1}{2},1])$ by Algorithm 4.3}
  \includegraphics[width=4cm]{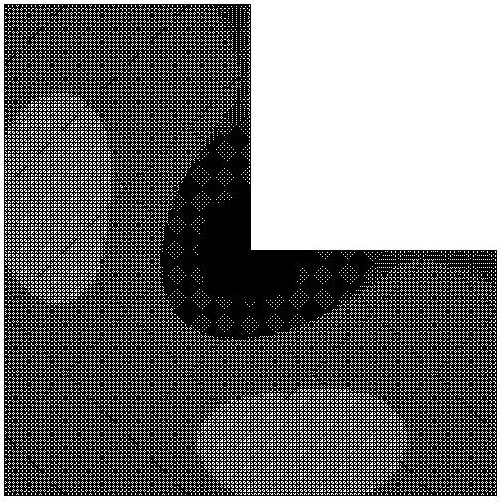}\\
  \caption{\small The adaptively refined mesh with 27537 degrees of freedom for $\lambda_2$ on $([0,1]\times [0,\frac{1}{2}])\bigcup ([0,\frac{1}{2}]\times
[\frac{1}{2},1])$ by Algorithm 4.3}
  \includegraphics[width=4cm]{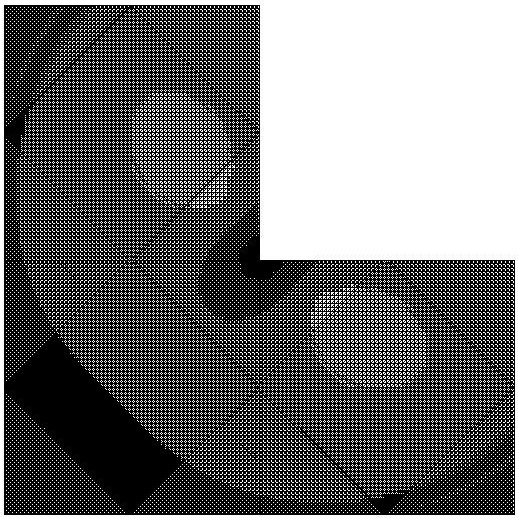}\\
  \caption{\small The adaptively refined mesh with 26341 degrees of freedom for $\lambda_3$ on $([0,1]\times [0,\frac{1}{2}])\bigcup ([0,\frac{1}{2}]\times
[\frac{1}{2},1])$ by Algorithm 4.3}
\end{figure}

\begin{figure}[htp]
  \centering
  \includegraphics[width=6cm]{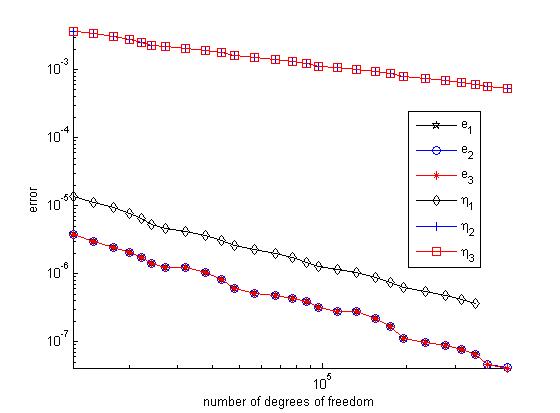}\\
  \caption{\small  The curves of error and the a posteriori estimators for the first eigenvalue of Example 4.2 in log-log scale}
\end{figure}

\begin{figure}[htp]
  \centering
  \includegraphics[width=6cm]{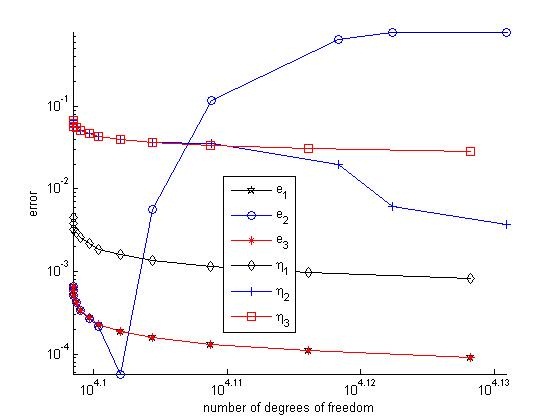}\\
  \caption{\small The curves of error and the a posteriori estimators for the second eigenvalue of Example 4.2 in log-log scale}
\end{figure}

\begin{figure}[htp]
   \centering
  \includegraphics[width=6cm]{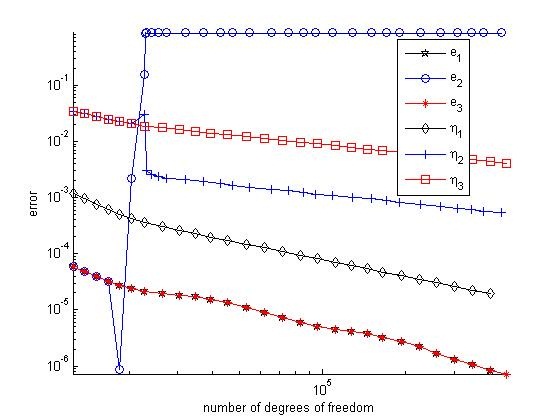}\\
  \caption{\small The curves of error and the a posteriori estimators for the third eigenvalue of Example 4.2 in log-log scale}
\end{figure}

Table 4 shows that the efficiency of Algorithm 4.3 is higher than that of Algorithm 4.1.
Table 5 illustrates that compared with Algorithm 4.2, our method can be used to solve approximations for any eigenvalue on the L-shaped domain
although the CPU time of Algorithm 4.2 is less. In detail, we list some results in Table 6
from which we can see that the approximate eigenvalues $\lambda_{2*}^{h_{l}}$ and
$\lambda_{3*}^{h_{l}}$ obtained by Algorithm 4.2 are not good
approximations for the second and the third exact eigenvalue but
are close to the first exact one (in Figs 8 and 9, the curves of $e_2$ also reflect this phenomenon),
and the reason we have analyzed in Example 4.1.

\noindent {\bf Remark 4.1.} In our computation, Algorithm 4.2 and Algorithm 4.3 are performed by solving an eigenvalue problem on the
coarsest mesh and then projecting the discrete eigenfunction on the subsequent finer meshes where only one linear system needs to be
solved. As for multiple eigenvalues, in recent publication \cite{dai2}, the authors considered the approximation of the eigenspace
and summed all the residuals up as for the a posteriori error estimate of eigenspace for the first type adaptive algorithm. Applying
this method to Algorithm4.2 and Algorithm 4.3 would be expected.

\section{Concluding Remarks}
\indent In this paper we propose and analyze an a posteriori error estimator
for the finite element multi-scale discretization approximation of the Steklov eigenvalue problem. Based on the a posteriori
error estimates, we design an adaptive algorithm of shifted-inverse iteration type. With this adaptive algorithm we can seek efficiently approximations of any eigenpair of the Steklov eigenvalue problem. \\
\indent The techniques used in this paper can also be applied to
finite elements approximation of second order
self-adjoint elliptic eigenvalue problems. The global upper and
lower bound of the error in \cite{li2} can also be obtained by using the argument
in this paper, moreover, the local lower bound of the error can also be
derived.


\vskip0.1in
\baselineskip 15pt
\renewcommand{\baselinestretch}{1.15}


\end{document}